\documentclass[11pt,a4paper]{amsart}
\usepackage{url}
\usepackage[T1]{fontenc}
\usepackage{amsmath,amsthm,amssymb}
\newcommand{\coloneqq}{:=}%
\theoremstyle{plain}
\newtheorem{thm}{Theorem}
\newtheorem{lem}[thm]{Lemma}
\newtheorem{defi}[thm]{Definition}
\newtheorem{prop}[thm]{Proposition}
\theoremstyle{remark}
\newtheorem{rem}{Remark}
\author{Fran{\c c}ois Martin}
\address{Universit\'e Blaise Pascal -- Clermont-Ferrand\\
Laboratoire de Math\'ematiques Pures\\
Les C\'ezeaux\\ F--63177 Aubi\`ere cedex\\ France}
\email{Francois.Martin@math.univ-bpclermont.fr}
\author{Emmanuel Royer}
\address{Universit\'e Blaise Pascal -- Clermont-Ferrand\\
Laboratoire de Math\'ematiques Pures\\
Les C\'ezeaux\\ F--63177 Aubi\`ere cedex\\ France}
\email{Emmanuel.Royer@math.univ-bpclermont.fr}
\title{Rankin-Cohen brackets on quasimodular forms}
\date{\today}
\thanks{Both authors are partially supported by the ANR project ``Modunombres'' and the BQR project ``Nomex''. The final version of this paper has been written at the \emph{Instituto Nacional de Matem\'atica Pura et Aplicada} in Rio de Janeiro. Both authors thank Hossein Movasati for the invitation and the institution for excellent working conditions.}
\newcommand{\C}{\mathbb{C}}
\DeclareMathOperator{\De}{D}
\newcommand*{\ic}{\mathrm{i}}
\newcommand{\qm}[2]{\widetilde{M}_{#1}^{\leq #2}}
\newcommand{\qme}[2]{\widetilde{M}_{#1}^{#2}}
\newcommand{\Modu}[1]{M_{#1}}
\newcommand{\Parab}[1]{S_{#1}}
\newcommand{\sldz}{\mathrm{SL}(2,\Z)}
\newcommand{\pk}{\mathcal{H}}
\newcommand{\trsf}[3]{(#1\underset{#2}{\vert}#3)}
\DeclareMathOperator{\Vect}{Vect}
\DeclareMathOperator{\X}{X}
\def\vec#1{\ensuremath{\mathchoice
           {\mbox{\boldmath$\displaystyle\mathbf{#1}$}}
           {\mbox{\boldmath$\textstyle\mathbf{#1}$}}
           {\mbox{\boldmath$\scriptstyle\mathbf{#1}$}}
           {\mbox{\boldmath$\scriptscriptstyle\mathbf{#1}$}}}}
\newcommand{\Z}{\mathbb{Z}}
\begin{document}
\begin{abstract}
We give the algebra of quasimodular forms a collection of Rankin-Cohen operators. These operators extend those defined by Cohen on modular forms and, as for modular forms, the first of them provides a Lie structure on quasimodular forms. They also satisfy a ``Leibniz rule'' for the usual derivation. Rankin-Cohen operators are useful for proving arithmetical identities. In particular, we explain why Chazy equation has the exact shape it has. %
\end{abstract}
\maketitle%
\section*{Introduction}
The purpose of this paper is to present a generalisation for quasimodular forms of the Rankin-Cohen brackets for modular forms: for each $n \geq 0$, $k, \ell, s, t$ positive integers, we define bilinear differential operators $[\,,\,]_n$ sending $\qm ks\times \qm {\ell}t$ to $\qm{k+\ell+2n}{s+t}$. We have denoted $\qm ks$ the vector space of quasimodular forms of weight $k$ and depth less or equal than $s$ on $\sldz$ (see section~\ref{sec:qmf} for the definitions). %

We give a quite precise description of the image of this bilinear form in terms of modular and parabolic forms. This allows us to obtain efficiently classical differential equations and arithmetical identities. %

Then we prove that the Rankin-Cohen brackets satisfy the ``Leibniz rule'' for the normalized usual derivation ($\De\coloneqq\frac{\mathrm{d}}{2\pi\ic\mathrm{d} z}$): $\De[f, g]_n = [\De f, g]_n + [f, \De g]_n$. %

The first section is a presentation of the definitions and classical results concerning quasimodular forms and Rankin-Cohen brackets on modular forms. %

In the second section, we prove the following theorem. %
\begin{thm}\label{thm:RCQM}%
Let $k,\ell$ in $\Z_{>0}$, $s\in\{0,\dotsc,\lfloor k/2\rfloor\}$, $t\in\{0,\dotsc,\lfloor \ell/2\rfloor\}$
and $n\in\Z_{\geq 0}$. Define %
\[%
\Phi_{n;k,s;\ell,t}(f,g)\coloneqq\sum_{r=0}^n(-1)^r\binom{k-s+n-1}{n-r}\binom{\ell-t+n-1}{r}\De^{r}f\De^{n-r}g. %
\]
Then %
\[%
\Phi_{n;k,s;\ell,t}\left(\qm{k}{s},\qm{\ell}{t}\right)\subset\qm{k+\ell+2n}{s+t}. %
\]
\end{thm}

In some case we get a more precise description in terms of the spaces of modular forms $\Modu{k}$ and the spaces of parabolic forms $\Parab{k}$. %
\begin{prop}\label{prop_struct}%
Under the hypothesis of theorem~\ref{thm:RCQM}, if $n>0$ then %
\[%
\Phi_{n;k,s;\ell,t}\left(\qm{k}{s},\qm{\ell}{t}\right)\in\Parab{k+\ell+2n}\oplus\bigoplus_{j=1}^{s+t}\De^j\Modu{k+\ell+2n-2j}. %
\]
If moreover $n>s+t$, then  %
\[%
\Phi_{n;k,s;\ell,t}\left(\qm{k}{s},\qm{\ell}{t}\right)\in\Parab{k+\ell+2n}\oplus\bigoplus_{j=1}^{s+t-1}\De^j\Modu{k+\ell+2n-2j}\oplus \De^{s+t}\Parab{k+\ell+2n-2s-2t}.
\]
The same conclusion holds if $n=s+t$ and $f$ or $g$ vanishes at infinity. %
\end{prop}

\begin{rem}%
This notion is consistent with the one for modular forms, the standard Rankin-Cohen bracket of $f \in M_k$ and $g \in M_{\ell}$ is $\Phi_{n;k,0;\ell,0}(f,g)$ (see paragraph \ref{par:RCmod}). %
\end{rem}
\begin{rem}\label{rem:int}%
For $n \geq 0$, a bilinear differential operator $\Psi$ sending $\qm ks\times \qm {\ell}t$ to $\bigcup_v\qm{k+\ell+2n}{v}$ is necessarily (for weight compatibility reasons) a linear combination of $(f, g) \mapsto \De^{r}f \De^{n-r}f$, $r \in \{0, \ldots,n\}$. Such a differential operator sends in principle $\qm ks \times \qm{\ell}t$ to $\qm{k+\ell+2n}{s+t+n}$ (see lemma \ref{lem:derqm}). So the operator $\Phi$ introduced before has the advantage of reducing the depth of the quasimodular form obtained, and it was not obvious that such an operator was existing. %
\end{rem}
\begin{rem}
 Theorem~\ref{thm:RCQM} is valid for quasimodular forms on any subgroup of finite index in $\sldz$. %
\end{rem}

In the third section, we show that the behaviour of this operator under derivation is natural. %
\begin{thm}\label{thm:deriv}
Under the hypothesis of theorem~\ref{thm:RCQM}, for all $f\in\qm{k}{s}$ and $g\in\qm{\ell}{t}$, %
\[%
\De\Phi_{n;k,s;\ell,t}(f,g)=\Phi_{n;k,s;\ell+2,t+1}(f,\De g)+\Phi_{n;k+2,s+1;\ell,t}(\De f,g). %
\]
\end{thm}
\begin{rem}
For $f$ of weight $k$ and \emph{exact depth} $s$ and $g$ of weight $\ell$ and \emph{exact depth} $t$, we write $[f,g]_n$ instead of $\Phi_{n;k,s;\ell,t}(f,g)$. Recall (see proposition \ref{pro:stabderqm}) that if $h$ has weight $w>0$ and depth $d$ then $\De h$ has weight $w+2$ and depth $d+1$. The following theorem may then be rewritten as %
\[%
\De[f,g]_n=[\De f,g]_n+[f,\De g]_n.%
\]
\end{rem}

For modular forms, Zagier, Cohen and Manin showed \cite{MR1418868} that the sum of Rankin-Cohen brackets defines an associative product on the algebra $M=\prod_{k \geq 0} M_k$. In a recent paper, Bieliavski, Tang and Yao \cite{MR2319770} showed that this sum is isomorphic to the standard Moyal product. Do the Rankin-Cohen brackets for quasimodular forms introduced here have such a geometric interpretation ? %

The existence of Rankin-Cohen brackets (thanks to proposition~\ref{prop_struct}) provides a new tool to obtain arithmetical identities. For example, we recover the Ramanujan differential equations, Chazy differential equation (and explain why such a differential equation has to exist), van der Pol equality and Niebur equality. As usual, define for $h\geq 2$ the Eisenstein series:
\begin{equation}\label{eq_eisenstein}%
E_h(z)\coloneqq 1-\frac{2h}{B_h}\sum_{n=1}^{+\infty}\sigma_h(n)\exp(2\pi\ic nz)%
\end{equation}
where $B_h$ is the Bernoulli number and %
\[%
 \sigma_h(n)\coloneqq\sum_{d\mid n}d^h. %
\]

One of the three Ramanujan equations is %
\[%
 \De E_2=-\frac{1}{12}(E_4-E_2^2). %
\]
It is a direct consequence of %
\[%
 [E_2,\Delta]_1=\Delta E_4
\]
where $\Delta$ is the unique primitive form of weight $12$ on $\sldz$. If we write $\tau(n)$ for the $n$th coefficient of $\Delta$, Niebur \cite{MR0382171} equality is %
\[%
 \tau(n)=n^4\sigma_1(n)-24\sum_{a=1}^{n-1}(35a^4-52a^3n+18a^2n^2)\sigma_1(a)\sigma_1(n-a) %
\]
and it follows from %
\[%
 [E_2,E_2]_4=-48\Delta.
\]
Van der Pol \cite{MR0042599} equality is %
\[%
 \tau(n)=n^2\sigma_3(n)+60\sum_{a=1}^{n-1}a(9a-5n)\sigma_3(a)\sigma_3(n-a). %
\]
It follows from %
\[%
 [E_4,\De E_4]_1=960\Delta. %
\]
Many examples of the two previous type are given in \cite{RaSa}. Finally, a quite astonishing equality is Chazy differential equation. Its usual form is %
\[%
\De^3E_2=E_2\De^2E_2-\frac{3}{2}(\De E_2)^2 %
\]
and it follows from %
\begin{equation}\label{eq_chazyrc}%
 \left[[K,\Delta]_1,\Delta\right]_1=24\Delta K^2 %
\end{equation}
where $K=[E_2,\Delta]_1$. The most outer bracket is on modular forms since it may be shown that $[K,\Delta]_1$ has depth $0$. That such a differential equation has to exist is a consequence of the following proposition that we prove using Rankin-Cohen brackets. %
\begin{prop}\label{prop_dern}
 Let $n\geq 0$ and $r\in\{0,\dotsc,n\}$. Then %
\[%
 \De^rE_2\De^{n-r}E_2\in\bigoplus_{\substack{j=0\\ j\equiv n\pmod{2}}}^{n-4}\De^j\Parab{2n+4-2j}\oplus\C\De^{n}E_{4}\oplus\C\De^{n+1}E_2. %
\]
In particular, $[E_2,E_2]_0\in\C E_4+\C\De E_2$, $[E_2,E_2]_2\in\C\De^2 E_4$, $[E_2,E_2]_4\in\C\Delta$ and %
\[%
 [E_2,E_2]_{2n}\in\Parab{4(n+1)}\oplus\De^2\Parab{4n} \quad\text{if $n\geq 3$.}
\]
\end{prop}
Indeed for $n=2$, this proposition implies that both quasimodular forms $E_2\De^2E_2$ and $(\De E_2)^2$ are in $\C\De^2E_{4}\oplus\C\De^3E_2$. Hence $\Vect(E_2\De^2E_2,(\De E_2)^2)=\Vect(\De^2E_{4},\De^3 E_2)$ and $\De^3 E_2$ is a linear combination of $E_2\De^2E_2$ and $(\De E_2)^2$: this is the shape of Chazy equation. %
\section{Overview}
\subsection{Quasimodular forms}\label{sec:qmf}
In this section, we introduce usual definitions and notations and recall some useful properties of quasimodular forms. For a more detailed introduction, see \cite[\S 17]{MR2186573}.

We introduce the following notations: let $\gamma = \begin{pmatrix} a & b \\ c & d \end{pmatrix} \in \sldz$ and $z \in \pk$, we define %
\[%
\X(\gamma,z)\coloneqq\frac{c}{cz+d} %
\]
and %
\[%
\X(\gamma) \colon z \mapsto \X(\gamma,z). %
\]
As usual, the complex upper half-plane is denoted by $\pk$. For $k\geq 0$, $f \colon \pk \rightarrow \C$ and $\gamma = \begin{pmatrix} a & b \\ c & d \end{pmatrix} \in \sldz$ the function $\trsf fk{\gamma}$ is defined by $\trsf fk{\gamma}(z) = (cz+d)^{-k} f(\gamma z)$.

\begin{defi}
Let $k\in\mathbb{Z}_{\geq 0}$ and $s\in\mathbb{Z}_{\geq 0}$. An holomorphic function $f \colon \pk \to \C $ is a quasimodular form of weight $k$, depth $s$ (over $\sldz$) if there exist holomorphic functions $Q_0(f)$, $Q_1(f)$, $\dotsc$, $Q_s(f)$ on $\pk$ such that %
\begin{equation}\label{eq_cqm}%
\trsf fk{\gamma} = \sum_{i=0}^s Q_i(f) \X(\gamma)^i %
\end{equation}
for all $\bigl(\begin{smallmatrix} a & b\\ c & d\end{smallmatrix}\bigr)\in\sldz$ and such that $Q_s(f)$ is not identically vanishing and $f$ has no negative terms in its Fourier expansion. By convention, the $0$ function is a quasimodular form of depth $-\infty$ and any weight.
\end{defi}
\begin{rem}
 Taking $\gamma=\bigl(\begin{smallmatrix} 1 & 0\\ 0 & 1\end{smallmatrix}\bigr)$ and $\gamma=\bigl(\begin{smallmatrix} 1 & 1\\ 0 & 1\end{smallmatrix}\bigr)$ in \eqref{eq_cqm} implies that $f$ is periodic of period $1$ hence has a Fourier expansion. The definition requires this Fourier expansion to be of the shape %
\[%
 f(z)=\sum_{n=0}^{+\infty}\widehat{f}(n)e^{2\pi\ic nz}.
\]
\end{rem}
The set of quasimodular forms of weight $k$ and depth $s$ is denoted by $\widetilde{M}_{k}^{s}$. It is often more convenient to use the $\C$-vectorial space of quasimodular forms of weight $k$ and depth less or equal than $s$, which is denoted by $\qm ks$. It can be shown that there are no quasimodular forms (except $0$) of negative weight or of depth $s>k/2$ \cite[lemme 120]{MR2186573}. Hence we extend our notation by defining $M_k^{\leq s}=\{0\}$ if $k<0$ and $M_k^{\leq s}=M_k^{\leq k/2}$ if $s>k/2$. %

\begin{rem}%
With this definition, the space $M_k$ of modular forms of weight $k$ for $\sldz$ is exactly the space $\qm k0$. %
\end{rem}
\begin{rem}%
A basic example of quasimodular form which is not a modular form is $E_2$ defined in \eqref{eq_eisenstein}. It satisfies for all $\gamma \in \sldz$ the transformation property \[\trsf{E_2}2{\gamma} = E_2 + \dfrac 6{\pi\ic} \X(\gamma),\] proving that $E_2 \in \qme 21$ (see e.g., \cite[lemme 19]{MR2186573}). %
\end{rem}

The space $\displaystyle{\qme *{} = \bigcup_{k, s} \qm k{s}}$ is equipped with a natural filtered-graded algebra structure (the graduation according to the weight, the filtration according to the depth). The canonical multiplication $(f, g) \longmapsto fg$ defines a morphism $\qm k{s}\times \qm {\ell}{t} \longrightarrow \qm{k+\ell}{s+t}$. %

If $f \in \qm ks$, the sequence $\left(Q_i(f)\right)_{i \in \Z}$ is defined by the quasimodularity condition \eqref{eq_cqm}, if $i \in \{0, \ldots, s\}$, and $Q_i(f)=0$ for $i \notin \{0, \ldots, s\}$. One can show that $Q_0(f)=f$ and $Q_i(f)\in\qm{k-2i}{s-i}$ \cite[Lemme 119]{MR2186573}. %

Quasimodular forms are the natural extension of modular forms into a stable by derivation space, because of the following proposition. %
\begin{prop}\label{pro:stabderqm}
If $k>0$, the normalized derivation $\De\coloneqq\dfrac{\mathrm{d}}{2\pi\ic\mathrm{d}z}$ maps $\qme ks$ to $\qme{k+2}{s+1}$. %
\end{prop}
For $r \in \Z_{\geq 0}$, write $f^{(r)}\coloneqq\De^r(f)$ and $f'=f^{(1)}$. The following lemma connects the transformation equation of $f$ and $f^{(r)}$. %
\begin{lem}\label{lem:derqm}%
Let $f\in\qm{k}{s}$. Then, %
\begin{equation}\label{for:derqm}%
\trsf{\De^{r}f}{k+2r}{\gamma}=\sum_{i=0}^{s+r}\left[\sum_{j=0}^r\frac{1}{(2\pi\ic)^j}j!\binom{r}{j}\binom{k+r-i+j-1}{j}\De^{r-j}Q_{i-j}(f)\right]\X(\gamma)^i %
\end{equation}
for all $r\in\Z_{\geq 0}$ and $\gamma\in\Gamma$. %
\end{lem}
\begin{proof} %
The result is obtained inductively on $r$: it is obvious for $r=0$, and for the induction suppose that for $r \geq 0$, formula~\eqref{for:derqm} holds. Let $g=f^{(r)}$. For $i \in \Z$ we have %
\begin{equation}\label{eq_veine}%
Q_i(g) =\sum_{j=0}^r\frac{1}{(2\pi\ic)^j}j!\binom{r}{j}\binom{k+r-i+j-1}{j}Q_{i-j}(f)^{(r-j)}\in\qm{k+2r-2i}{s+r-i}. %
\end{equation}
Then using proposition~\ref{pro:stabderqm} (which implies that $f^{(r+1)}\in \qm{k+2r+2}{r+s+1}$) and lemma 118 of \cite{MR2186573} we find %
\[%
\trsf{f^{(r+1)}}{k+2r+2}{\gamma}=\sum_{i=0}^{s+r+1} \left(Q_i(g)' +\frac{k+2r-i+1}{2\pi\ic}Q_{i-1}(g)\right)\X(\gamma)^i. %
\]
From~\eqref{eq_veine} we compute %
\begin{multline*}%
Q_i(g)'+\frac{k+2r-i+1}{2\pi\ic}Q_{i-1}(g)=\\%
Q_i(f)^{(r+1)}+\frac{k+2r-i+1}{(2\pi\ic)^{r+1}}r!\binom{k+2r-i}{r}Q_{i-r-1}(f) + \sum_{j=1}^r\frac{1}{(2\pi\ic)^j}Q_{i-j}(f)^{(r+1-j)} \\%
\times\left(\frac{r!}{(r-j)!}\binom{k+r-i+j-1}{j}+\frac{(k+2r-i+1)r!}{(r+1-j)!}\binom{k+r-i+j-1}{j-1}\right). %
\end{multline*}
Formula~\eqref{for:derqm} for $r+1$ instead of $r$ follows by expanding the binomial coefficients. %
\end{proof}

Finally, we shall need the following structure result. For completness, we provide a short proof that should convince that the theory requires $E_2$. %
\begin{prop}\label{prop_struqm}%
 Quasimodular forms can be expressed as linear combinations of derivatives of modular forms and $E_2$ : %
\[%
 \qm{k}{k/2}=\bigoplus_{i=0}^{k/2-1}\De^i\Modu{k-2i}\oplus\C\De^{k/2-1}E_2. %
\]
\end{prop}
\begin{proof}%
 We proceed by descent on the depth. If $f$ has weight $k$ and depth $s$, we would like to have a modular form $g$ such that $f-\De^sg$ has depth strictly less than $s$.  For any $g\in\Modu{k-2s}$, multiple use of differentiation theorem \cite[Lemme 118]{MR2186573} lead to %
\begin{equation}\label{eq_qsds}%
 Q_s(D^sg)=\left(\frac{1}{2\pi\ic}\right)^ss!\binom{k-s-1}{s}g. %
\end{equation}
If $\binom{k-s-1}{s}\neq 0$, which happens if $s<k/2$, we can choose %
\[%
 g=(2\pi\ic)^s\frac{(k-2s-1)!}{(k-s-1)!}Q_s(f)\in\Modu{k-2s}. %
\]
For $s=k/2$, we use %
\[%
 Q_{k/2}(D^{k/2-1}E_2)=\left(\frac{1}{2\pi\ic}\right)^{k/2-1}\left(\frac{k}{2}-1\right)!\frac{6}{\pi\ic} %
\]
and choose %
\[%
 \alpha=\frac{\pi\ic}{6}\cdot\frac{(2\pi\ic)^{k/2-1}}{\left(\frac{k}{2}-1\right)!}Q_{k/2}(f)\in\Modu{0}=\C %
\]
to obtain %
\[%
 f-\alpha\De^{k/2}E_2\in\qm{k}{k/2-1}. %
\]
\end{proof}

\subsection{Usual Rankin-Cohen brackets for modular forms}\label{par:RCmod}
The Rankin-Cohen brackets have been introduced by Cohen after a work of Rankin. These are bilinear differential operators, whose main property is to preserve modular forms. More precisely, let $\Gamma$ be a finite index subgroup of $\sldz$. We write $M_k(\Gamma)$ for the space of modular forms of weight $k$ over $\Gamma$. For each $n \geq 0$, $(f, g) \in M_k(\Gamma)\times M_{\ell}(\Gamma)$, define the $n$-Rankin-Cohen bracket of $f$ and $g$ by %
\begin{equation}%
[f, g]_n = \sum_{r=0}^n (-1)^r \binom{k+n-1}{n-r}\binom{\ell+n-1}{r}\De^{r}f\De^{n-r}g.%
\end{equation}
Then $[f, g]_n \in M_{k+\ell+2n}(\Gamma)$. Moreover, if $\Phi$ is a bilinear differential operator sending $M_k(\Gamma)\times M_{\ell}(\Gamma)$ to $M_{k+\ell+2n}(\Gamma)$ for all $\Gamma\subset \sldz$ a finite index subgroup, then (up to constant) $\Phi(f,g)=[f, g]_n$. For an overview of Rankin-Cohen brackets including a proof of these results\footnote{The uniqueness result needs explanations: it is proved by using only algebraic arguments, the demonstration does not depend on the group $\Gamma$ or on growth conditions. Of course, it is possible that for some fixed group $\Gamma$ the uniqueness result does not hold (for instance if $M_k(\Gamma) = \{ 0 \}$ !).}, see for instance \cite{Zag94}, \cite{Zag92} or \cite{MR2186573}. %

Rankin-Cohen brackets appear to be useful in various mathematical domains as for instance invariant theory (\cite{Unt96} and \cite{Cho01}) or non-commutative geometry \cite{Yao}. %

\section{Rankin-Cohen brackets}

We prove our main result (theorem~\ref{thm:RCQM}). For $n \geq 0$ and any sequence $\vec{a}=(a_r)_{0\leq r\leq n}$, the bilinear forms we study take the form %
\[%
\Phi_{\vec{a}}(f,g)=\sum_{r=0}^na_r\De^rf\De^{n-r}g. %
\]
We first establish a sufficient condition on $\vec{a}$ (lemma~\ref{lem:fico}). For $s$, $t$ and $n$ nonnegative integers, we introduce the set %
\[%
\mathcal{E}(s,t,n)=\bigl\{(u,v,\alpha,\beta)\in\Z_{\geq 0}^4\colon u\leq s,\,v\leq t,\,\alpha+\beta\leq u+v+n-s-t-1\bigr\}. %
\]

\begin{lem}\label{lem:fico}%
Let $k,\ell$ in $\Z_{>0}$, $s\in\{0,\dotsc,\lfloor k/2\rfloor\}$, $t\in\{0,\dotsc,\lfloor \ell/2\rfloor\}$ and $n\in\Z_{>0}$. For $\vec{a}=(a_r)_{0\leq r\leq n}$ satisfying %
\[%
\sum_{r=0}^na_r\binom{r}{\alpha}\binom{n-r}{\beta}(k+r-u-1)!(\ell+n-r-v-1)!=0%
\]
for all $(u,v,\alpha,\beta)\in\mathcal{E}(s,t,n)$, one has %
\[%
\Phi_{\vec{a}}\left(\qm{k}{s},\qm{\ell}{t}\right)\subset\qm{k+\ell+2n}{s+t}. %
\]
\end{lem}
\begin{proof}%
Let $f\in\qm{k}{s}$ and $g\in\qm{\ell}{t}$. From lemma~\ref{lem:derqm} we deduce %
\begin{align*}
\trsf{\Phi_{\vec{a}}(f,g)}{k+\ell+2n}{\gamma}&=\sum_{r=0}^{n}a_r\trsf{f^{(r)}}{k+2r}{\gamma}\trsf{g^{(n-r)}}{\ell+2(n-r)}{\gamma}\\
&=\sum_{i=0}^{s+t+n}C(\vec{a};i)(f,g)\X(\gamma)^i %
\end{align*}
with %
\begin{multline}\label{eq_defveca}%
C(\vec{a};i)(f,g)=\sum_{\substack{(i_1,i_2)\in\Z_{\geq 0}^2\\ i_1+i_2=i}}\sum_{r=0}^na_r\sum_{j_1=0}^r\left(\frac{1}{2\pi\ic}\right)^{j_1}j_1!\binom{r}{j_1}\binom{k+r-i_1+j_1-1}{j_1}\\\times\sum_{j_2=0}^{n-r}\left(\frac{1}{2\pi\ic}\right)^{j_2}j_2!\binom{n-r}{j_2}\binom{\ell+n-r-i_2+j_2-1}{j_2}Q_{i_1-j_1}(f)^{(r-j_1)}Q_{i_2-j_2}(g)^{(n-r-j_2)}. %
\end{multline}
It follows that $\Phi_{\vec{a}}(f,g)\in\qm{k+\ell+2n}{s+t}$ if and only if $C(\vec{a};s+t+i)=0$ for all $i\in\{1,\dotsc,n\}$. This is easily seen to be equivalent to %
\begin{multline*}%
\sum_u\sum_v\sum_{\substack{(\alpha,\beta)\in\Z_{\geq 0}^2\\\alpha+\beta=n+u+v-s-t-i}}\left(\frac{1}{2\pi\ic}\right)^{n-\alpha-\beta}\sum_ra_r(r-\alpha)!(n-r-\beta)!\\\times\binom{r}{\alpha}\binom{n-r}{\beta}\binom{k+r-u-1}{r-\alpha}\binom{\ell+n-r-v-1}{n-r-\beta}Q_u(f)^{(\alpha)}Q_v(g)^{(\beta)}=0 %
\end{multline*}
for all $i\in\{1,\dotsc,n\}$, the sets of summation being determined by the binomial coefficients. Hence, $\Phi_{\vec{a}}\left(\qm{k}{s},\qm{\ell}{t}\right)\subset\qm{k+\ell+2n}{s+t}$ is implied by %
\begin{equation}\label{eq_etoile}%
\sum_ra_r\binom{r}{\alpha}\binom{n-r}{\beta}(k+r-u-1)!(\ell+n-r-v-1)!=0 %
\end{equation}
for all $(u,v,\alpha,\beta)\in\mathcal{E}(s,t,n)$. %
\end{proof}

\begin{rem}\label{rem:uni}
The statement of the previous lemma is in fact an equivalence, if we ask $\Phi_{\vec{a}}$ to satisfy $\Phi_{\vec{a}}\left(\qm{k}{s}(\Gamma),\qm{\ell}{t}(\Gamma)\right)\subset\qm{k+\ell+2n}{s+t}(\Gamma)$ for each finite index subgroup $\Gamma$ of $\sldz$: indeed for $\{a(u, v,\alpha, \beta)\}$ a non identically zero family of complex numbers, if $\Psi \colon (f,g) \mapsto\sum_{(u,v,\alpha,\beta)\in\mathcal{E}(s,t,n)} a(u,v,\alpha, \beta)Q_u(f)^{(\alpha)}Q_v(g)^{(\beta)}$ satisfy $\Psi(\qm ks(\Gamma), \qm{\ell}t(\Gamma)) = 0$, then exists $M>0$ such that the minimum of $\dim(\qm ks(\Gamma))$ and  $\dim(\qm{\ell}t(\Gamma))$ is strictly smaller than $M$. However, as for modular forms, for each $A>0$ exists $\Gamma$ a finite index subgroup of $\sldz$ such that $\dim \qm ks(\Gamma) > A$ and $\dim \qm {\ell}t(\Gamma) > A$ (recall that $k, \ell \in \Z_{>0}$). %
\end{rem}

We shall now give a necessary condition for $\vec{a}$ satisfying the condition of lemma~\ref{lem:fico}. %
\begin{lem}\label{lem:uniq}%
Let $k,\ell$ in $\Z_{>0}$, $s\in\{0,\dotsc,\lfloor k/2\rfloor\}$, $t\in\{0,\dotsc,\lfloor \ell/2\rfloor\}$ and $n\in\Z_{>0}$. If $\vec{a}=(a_r)_{0\leq r\leq n}$ satisfies %
\[%
\sum_{r=0}^na_r\binom{r}{\alpha}\binom{n-r}{\beta}(k+r-u-1)!(\ell+n-r-v-1)!=0 %
\]
for all $(u,v,\alpha,\beta)\in\mathcal{E}(s,t,n)$, then there exists $\lambda\in\C$ such that %
\[%
a_r=\lambda(-1)^r\binom{k+n-s-1}{n-r}\binom{\ell+n-t-1}{r} %
\]
for all $r\in\{0,\dotsc,n\}$. %
\end{lem}
\begin{proof}
Define $\vec{b}=(b_r)_{0\leq r\leq n}$ by %
\[%
b_r=a_r(k+r-s-1)!(\ell+n-r-t-1)! %
\]
for all $r$. Then %
\[%
\sum_{r=0}^nb_r\binom{r}{\alpha}\binom{n-r}{\beta}\binom{k+r-u-1}{s-u}\binom{\ell+n-r-v-1}{t-v}=0 %
\]
for all $(u,v,\alpha,\beta)\in\mathcal{E}(s,t,n)$. Choosing $u=s$, $t=v$ and $\beta=0$ leads to $F^{(\alpha)}(1)=0$ for all $\alpha\in\{0,\dotsc,n-1\}$ where $F$ is the generating (polynomial) function of $\vec{b}$ defined by %
\[%
F(x)=\sum_{r=0}^nb_rx^r. %
\]
This implies the existence of $\mu\in\C$ such that $F(x)=\mu(x-1)^n$ and thus $b_r=\mu(-1)^r\binom{n}{r}$. The result follows by defining %
\[%
\lambda=\mu\frac{n!}{(k-s+n-1)!(\ell-t+n-1)!}. %
\]
\end{proof}
We obtain the existence of the Rankin-Cohen operator for quasimodular forms in showing that the vector $\vec{a}$ we found in lemma~\ref{lem:uniq} is admissible. %
\begin{lem}\label{lem:exist}%
Let $k,\ell$ in $\Z_{>0}$, $s\in\{0,\dotsc,\lfloor k/2\rfloor\}$, $t\in\{0,\dotsc,\lfloor \ell/2\rfloor\}$ and $n\in\Z_{>0}$. Let $\vec{a}=(a_r)_{1\leq r\leq n}$ be defined by %
\[%
a_r=(-1)^r\binom{k-s+n-1}{n-r}\binom{\ell-t+n-1}{r}. %
\]
Then %
\[%
\Phi_{\vec{a}}\left(\qm{k}{s},\qm{\ell}{t}\right)\subset\qm{k+\ell+2n}{s+t}. %
\]
\end{lem}
\begin{proof}
By lemma~\ref{lem:fico} it suffices to check that %
\begin{equation}\label{eq_tbp}%
\sum_{\substack{(r_1,r_2)\in\Z_{\geq 0}\times\Z_{\geq 0}\\ r_1+r_2=n}}\frac{(-1)^{r_1}}{r_1!r_2!}\binom{r_1}{\alpha}\binom{r_2}{\beta}\binom{k-u-1+r_1}{s-u}\binom{\ell-v-1+r_2}{t-v}=0 %
\end{equation}
for all $(u,v,\alpha,\beta)\in\mathcal{E}(s,t,n)$. Fix  $(u,v,\alpha,\beta)\in\mathcal{E}(s,t,n)$, then \eqref{eq_tbp} is the coefficient of order $n$ in the product $P_1(X)P_2(X)$ where %
\begin{align*}%
P_1(X) &= \sum_{r_1=0}^{+\infty}\frac{(-1)^{r_1}}{r_1!}\binom{r_1}{\alpha}\binom{k-u-1+r_1}{s-u}X^{r_1}\\%
P_2(X) &= \sum_{r_2=0}^{+\infty}\frac{1}{r_2!}\binom{r_2}{\beta}\binom{\ell-v-1+r_2}{t-v}X^{r_2}. %
\end{align*}
We have %
\[%
P_1(X)=\frac{X^\alpha}{\alpha!}Q_1^{(\alpha)}(X) %
\]
with %
\[%
Q_1(X)=\sum_{r_1=0}^{+\infty}\frac{(-1)^{r_1}}{r_1!}\binom{k-u-1+r_1}{s-u}X^{r_1} %
\]
and %
\[%
Q_1(X)=\frac{X^{-k+s+1}}{(s-u)!}R_1^{(s-u)}(X) %
\]
with %
\begin{align*}%
R_1(X)&=\sum_{r_1=0}^{+\infty}\frac{(-1)^{r_1}}{r_1!}X^{r_1+k-u-1}\\%
      &=X^{k-u-1}e^{-X}. %
\end{align*} %
We therefore may write $P_1(X)=\Pi_1(X)e^{-X}$ where $\Pi_1$ is a polynomial of degree $\alpha+s-u$. Similary, $P_2(X)=\Pi_2(X)e^{X}$ where $\Pi_2$ is a polynomial of degree $\beta+t-v$. It follows that $P_1P_2$ is a polynomial of degree $\alpha+\beta+s+t-u-v$. Finally, since, by definition, $\alpha+\beta-u-v<n-s-t$ we get \eqref{eq_tbp}. %
\end{proof}
\begin{rem}
With the help of the hypergeometric methods \cite[Chapter 3]{MR1379802}, we obtain that %
\[%
 \Pi_1(X)=(-1)^\alpha\sum_{r=\alpha}^{s-u+\alpha}\binom{k+\alpha-u-1}{k+r-s-1}\binom{r}{\alpha}\frac{X^{r}}{r!} %
\]
and %
\[%
 \Pi_2(X)=(-1)^\beta\sum_{r=\beta}^{t-v+\beta}(-1)^r\binom{\ell+\beta-v-1}{\ell+r-t-1}\binom{r}{\beta}\frac{X^{r}}{r!} %
\]
\end{rem}
Previous lemmas prove theorem \ref{thm:RCQM}. %
\section{Rankin-Cohen brackets and derivation}
In this section, we prove theorem~\ref{thm:deriv}. First, we remark that %
\begin{multline}\label{eq_1812} %
\Phi_{n;k,s;\ell,t}(f,g)'=\sum_{r=0}^{n-1}(-1)^r\Biggl[\binom{k-s+n-1}{n-r}\binom{\ell-t+n-1}{r}\\-\binom{k-s+n-1}{n-r-1}\binom{\ell-t+n-1}{r+1}\Biggr]f^{(r+1)}g^{(n-r)}\\+\binom{k-s+n-1}{n}fg^{(n+1)}+(-1)^{n}\binom{\ell-t+n-1}{n}f^{(n+1)}g. %
\end{multline}
Next, %
\begin{multline*}%
\Phi_{n;k,s;\ell+2,t+1}(f,g')=\binom{k-s+n-1}{n}fg^{(n+1)}\\-\sum_{r=0}^{n-1}(-1)^r\binom{k-s+n-1}{n-r-1}\binom{\ell-t+n}{r+1}f^{(r+1)}g^{(n-r)}
\end{multline*}
so that %
\begin{multline}\label{eq_1813}
\Phi_{n;k+2,s+1;\ell,t}(f',g)+\Phi_{n;k,s;\ell+2,t+1}(f,g')=\\\binom{k-s+n-1}{n}fg^{(n+1)}+(-1)^n\binom{\ell-t+n-1}{n}f^{(n+1)}g\\+\sum_{r=0}^{n-1}(-1)^r\Biggl[\binom{k-s+n}{n-r}\binom{\ell-t+n-1}{r}\\-\binom{k-s+n-1}{n-r-1}\binom{\ell-t+n}{r+1}\Biggr]f^{(r+1)}g^{(n-r)} %
\end{multline}
and equality from \eqref{eq_1812} and \eqref{eq_1813} follows by expanding the binomial coefficients. %
\section{A more precise structure result}%
In this section, we prove proposition~\ref{prop_struct}. Let $n>0$. If $f\in\qme{k}{s}$ and $g\in\qme{\ell}{t}$ then $\Phi_{n;k,s;\ell,t}(f,g)$ has weight $k+\ell+2n$ and depth less than $s+t$. Since $n>0$ this depth is not maximal since %
\[%
 s+t\leq\frac{k}{2}+\frac{\ell}{2}<\frac{k+\ell+2n}{2}. %
\]
Then it follows from proposition~\ref{prop_struqm} that %
\[%
 \Phi_{n;k,s;\ell,t}(f,g)\in\Modu{k+\ell+2n}\oplus\bigoplus_{j=1}^{s+t}\De^j\Modu{k+\ell+2n-2j}. %
\]
However, the definition of $\Phi_{n;k,s;\ell,t}(f,g)$ implies that its Fourier coefficient at $0$ is $0$ and since this is also true for derivatives of modular forms we get %
\[%
 \Phi_{n;k,s;\ell,t}(f,g)\in\Parab{k+\ell+2n}\oplus\bigoplus_{j=1}^{s+t}\De^j\Modu{k+\ell+2n-2j}. %
\]
The contribution to $\Phi_{n;k,s;\ell,t}(f,g)$ coming from %
\[%
 \Parab{k+\ell+2n}\oplus\bigoplus_{j=1}^{s+t-1}\De^j\Modu{k+\ell+2n-2j} %
\]
has depth strictly less than $s+t$. Hence %
\[%
 Q_{s+t}\left(\Phi_{n;k,s;\ell,t}(f,g)\right)=Q_{s+t}(\De^{s+t}g) %
\]
where $g\in\Modu{k+\ell+2n-2s-2t}$. Since %
\[%
Q_{s+t}(\De^{s+t}g)=(2\pi\ic)^{-s-t}\frac{(k+\ell+2n-s-t-1)!}{(k+\ell+2n-2s-2t-1)!}g %
\]
(see~\eqref{eq_qsds}), to prove that $g$ is parabolic we shall prove that the Fourier coefficient at $0$ of $Q_{s+t}\left(\Phi_{n;k,s;\ell,t}(f,g)\right)$ is $0$. From~\eqref{eq_defveca} we get %
\begin{multline}\label{eq_retourpuy}%
 Q_{s+t}\left(\Phi_{n;k,s;\ell,t}(f,g)\right)=\sum_u\sum_v\sum_{\substack{(\alpha,\beta)\in\Z_{\geq0}^2\\\alpha+\beta=n+u+v-s-t}}\left(\frac{1}{2\pi\ic}\right)^{n-\alpha-\beta}\sum_ra_r(r-\alpha)!(n-r-\beta)!\\\times\binom{r}{\alpha}\binom{n-r}{\beta}\binom{k+r-u-1}{r-\alpha}\binom{\ell+n-r-v-1}{n-r-\beta}Q_u(f)^{(\alpha)}Q_v(g)^{(\beta)}. %
\end{multline}
Since derivatives of quasimodular forms have Fourier coefficients vanishing at $0$, the only contribution to the Fourier coefficient of $Q_{s+t}\left(\Phi_{n;k,s;\ell,t}(f,g)\right)$ at $0$ is given by $(\alpha,\beta)=(0,0)$ in~\eqref{eq_retourpuy}. However, the summation is on $(\alpha,\beta)$ such that $\alpha+\beta=n+u+v-s-t$ and we have $n+u+v-s-t>0$ if $n>s+t$. %
 Thanks to \eqref{eq_retourpuy} we also see that if $f\in\qm{k}{s}$ and $g\in\qm{\ell}{t}$ satisfies $s+t>0$ and $\widehat{g}(0)=0$ then %
\[%
 \Phi_{s+t;k,s;\ell,t}(f,g)\in\Parab{k+\ell+2s+2t}\oplus\bigoplus_{j=1}^{s+t-1}\De^j\Modu{k+\ell+2s+2t-2j}\oplus\De^{s+t}\Parab{k+\ell}.
\]
\section{Applications}
An easy but useful consequence of the fact that $\De\Delta=\Delta E_2$ is the following lemma. %
\begin{lem}
Let $n\geq 0$. Let $f\in\qm{k}{s}$ and $g\in\qm{\ell}{t}$. There exists $h\in\qm{k+\ell+2n}{s+t}$ such that %
\[%
 \Phi_{n;k,s;\ell,t}(f,\Delta g)=\Delta h. %
\]
\end{lem}
For example, we have %
\[%
 \Phi_{1;k+12,s;12,0}(\Delta f,\Delta)=\Delta\Phi_{1;k,s;12,0}(f,\Delta).
\]

\subsection{Homogoneous products of derivatives of $E_2$}

In this section we prove proposition~\ref{prop_dern} by recursion on $n$. For $n=0$ we have $E_2^2=E_4+12\De E_2\in\C E_4\oplus\C\De E_2$. Assume that: %
\[
 \De^r E_2\De^{n-r}E_2\in\bigoplus_{\substack{j=0\\ j\equiv n\pmod{2}}}^{n-4}\De^j\Parab{2n+4-2j}\oplus\C\De^n E_4\oplus\C\De^{n+1}E_2 \quad(0\leq r\leq n).
\]

Deal first with the case where $n=2m$ is even. By recursion hypothesis, we have %
\begin{align*}%
 \De\left(\De^r E_2\De^{n-r}E_2\right)&=\De^r E_2\De^{n+1-r}E_2+\De^{r+1} E_2\De^{n-r}E_2\\ &\in\bigoplus_{\substack{j=0\\ j\equiv n\pmod{2}}}^{n-4}\De^{j+1}\Parab{2n+4-2j}\oplus\C\De^{n+1}E_4\oplus\C\De^{n+2}E_2.
\end{align*}
The set $\{\De^rE_2\De^{n-r}E_2,\, 0\leq r\leq n\}$ has $m+1$ distinct terms (corresponding to $0\leq r\leq m$). The set $\{\De^rE_2\De^{n+1-r}E_2,\, 0\leq r\leq n+1\}$ has also $m+1$ distinct terms (corresponding to $0\leq r\leq m$). It follows that %
\[%
 \left\{\De^r E_2\De^{n+1-r}E_2+\De^{r+1} E_2\De^{n-r}E_2,\, r\in\{0,\dotsc,m\}\right\}%
\]
and%
\[%
\left\{\De^r E_2\De^{n+1-r}E_2,\, r\in\{0,\dotsc,m\}\right\}%
\]
are basis of the same space with change of basis matrix given by %
\[%
 \begin{pmatrix}
1      & 0       & \dots   & \dots  & 0 \\
1      & 1       & \ddots  &        & \vdots  \\
0      & 1       & \ddots  & \ddots & \vdots  \\
\vdots & \ddots  & \ddots  & 1      & 0 \\
0      & \dots   & 0       & 1      & 2%
 \end{pmatrix}.
\]
It follows that for any $r\in\{0,\dotsc,m\}$ (hence any $r\in\{0,\dotsc,n\}$) we have %
\[%
 \De^r E_2\De^{n+1-r}E_2\in\bigoplus_{\substack{j=0\\ j\equiv n+1\pmod{2}}}^{n-3}\De^j\Parab{2n+6-2j}\oplus\C\De^{n+1} E_4\oplus\C\De^{n+2}E_2.
\]
We now deal with the case where $n=2m-1$ is odd. Again, by recursion hypothesis, we have %
\begin{align*}%
 \De\left(\De^r E_2\De^{n-r}E_2\right)&=\De^r E_2\De^{n+1-r}E_2+\De^{r+1} E_2\De^{n-r}E_2\\ &\in\bigoplus_{\substack{j=0\\ j\equiv n\pmod{2}}}^{n-4}\De^{j+1}\Parab{2n+4-2j}\oplus\C\De^{n+1}E_4\oplus\C\De^{n+2}E_2.
\end{align*}
The subspace generated by all the quasimodular forms $\De^{r} E_2\De^{n+1-r} E_2+\De^{r+1} E_2\De^{n-r}E_2$ when $r$ runs over $\{0,\dotsc,2m-1\}$ is the hyperplane %
\[%
 \left\{\sum_{r=0}^{2m}\alpha_r\De^r E_2\De^{2m-r} E_2 \vert \sum_{r=0}^{2m}(-1)^r\alpha_r=0\right\}
\]
hence it is sufficient for the proof of our recursion step to find a linear combination %
\[%
 \sum_{r=0}^{2m}\alpha_r\De^r E_2\De^{2m-r} E_2\in\bigoplus_{\substack{j=0\\ \text{$j$ even}}}^{2m-4}\De^{j}\Parab{4m+4-2j}\oplus\C\De^{2m}E_4\oplus\C\De^{2m+1}E_2 %
\]
with %
\[%
 \sum_{r=0}^{2m}(-1)^r\alpha_r\neq 0. %
\]
This is the step where we use Rankin-Cohen brackets. Since $[E_2,E_2]_{2m+2}\in\qm{4m+8}{2}$ we have $Q_2\left([E_2,E_2]_{2m+2}\right)\in\Parab{4m+4}$ (see~\eqref{eq_retourpuy} for the cuspidality). Equation~\eqref{eq_defveca} combined with the fact that $Q_1(E_2)$ is constant implies that %
\begin{multline}\label{eq_QdRC}
 Q_2\left([E_2,E_2]_{2m+2}\right)=\frac{24}{(2\pi\ic)^2}(2m+2)\De^{2m+1}E_2+\frac{4}{(2\pi\ic)^2}\Biggl[\\\sum_{r=2}^{2m+2}(-1)^r\binom{2m+2}{r}^2\binom{r}{2}\binom{r+1}{2}\De^{r-2}E_2\De^{2m+2-r}E_2\\+\sum_{r=1}^{2m+1}(-1)^r\binom{2m+2}{r}^2\binom{r+1}{2}\binom{2m+3-r}{2}\De^{r-1}E_2\De^{2m+1-r}E_2\Biggr]. %
\end{multline}
Let %
\[%
 \alpha_r(N)=2(-1)^r\binom{r}{2}\binom{N}{r}\binom{N}{r-1}(N+1-2r). %
\]
Equation~\eqref{eq_QdRC} gives %
\begin{multline*}%
 \sum_{r=2}^{2m+2}\alpha_r(2m+2)\De^{r-2}E_2\De^{2m+2-r}E_2 = \\%
 (2\pi\ic)^2Q_2\left([E_2,E_2]_{2m+2}\right)-24(2m+2)\De^{2m+1}E_2 \\%
\in \Parab{4m+4}\oplus\C\De^{2m+1}E_2. %
\end{multline*}
Let $\beta_r(N)=(-1)^r\alpha_r(N)$. We prove that%
\[%
 A(N)=\sum_{r=2}^N(-1)^r\alpha_r(N)=\sum_{r\in\mathbb{Z}}\beta_r(N) %
\]
is strictly negative (hence differs from $0$). Zeilberger's algorithm (\textit{e.g.}, on the open-source computer algebra system \textsf{Maxima}) \cite[Chapter 6]{MR1379802} provides a function $K(N,r)$ such that\footnote{Note that no algorithm is needed to check that $K(N,r)$ as defined in~\eqref{eq_defG} works.} %
\begin{multline*}%
 2(N+1)(2N-1)\beta_r(N)-N(N-1)\beta_r(N+1)=\\ %
K(N,r+1)\beta_{r+1}(N)-K(N,r)\beta_r(N). %
\end{multline*}
More precisely %
\begin{multline}\label{eq_defG}%
 K(N,r)=\\\frac{(r-2)(r-1)(N+1)[3N^3+8N^2(1-r)+N(4r^2-6r+3)-2r^2+4r-2]}{(N-2r+1)(N-r+1)(N-r+2)(N-1)}. %
\end{multline}
We deduce the recursive formula %
\[%
 \frac{A(N+1)}{A(N)}=\frac{2(N+1)(2N-1)}{N(N-1)} %
\]
which, since $A(2)=4$, implies %
\[%
 A(N)=-N(N-1)\binom{2N-2}{N-1}<0.
\]
Finally, we have found a function which belongs to the hyperplane. This completes the proof. %

\subsection{Niebur formula}
From proposition~\ref{prop_dern} we obtain %
\[%
 \Phi_{4;2,1;2,1}(E_2,E_2)\in\Parab{12}=\C\Delta. %
\]
The computation of the first coefficients gives $\Phi_{4;2,1;2,1}(E_2,E_2)=-48\Delta$. This is the differential equation proved by Niebur in~\cite{MR0382171} : %
\[%
 2^3\cdot3\Delta=18(\De^2E_2)^2+E_2\De^4E_2-16\De E_2\De^3E_2 %
\]
and comparing the Fourier expansions gives Niebur formula. %
\subsection{van der Pol formula}
From proposition~\ref{prop_struct} we obtain %
\[%
 \Phi_{1;4,0;6,1}(E_4,\De E_4)\in\Parab{12}. %
\]
The computation of the first coefficient gives $\Phi_{1;4,0;6,1}(E_4,\De E_4)=960\Delta$. This is the differential equation proved by van der Pol: %
\[%
 4E_4\De^2E_4-5(\De E_4)^2=960\Delta. %
\]
It leads to  %
\begin{align*}%
 \tau(n)&=n^2\sigma_3(n)+60\sum_{a+b=n}(4b-5a)b\sigma_3(a)\sigma_3(b)\\
        &=n^2\sigma_3(n)+60\sum_{a=1}^{n-1}(9a^2-13an+4n^2)\sigma_3(a)\sigma_3(n-a)\\
        &=n^2\sigma_3(n)+60\sum_{b=1}^{n-1}(9b^2-5bn)\sigma_3(a)\sigma_3(n-a)%
\end{align*}
and the summation of the two last equalities implies the van der Pol formula in its original form \cite[eq. (53)]{MR0042599}: %
\[%
 \tau(n)=n^2\sigma_3(n)+60\sum_{a=1}^{n-1}(2n-3a)(n-3a)\sigma_3(a)\sigma_3(n-a).
\]

\subsection{Chazy equation}
Recall that we proved at the end of the introduction that an equation of the shape
\[%
 \alpha E_2\De^2 E_2+\beta(\De E_2)^2=\De^3E_2 %
\]
has to exist. Coefficients $\alpha$ and $\beta$ can be computed by identifications of the first Fourier coefficients. Our aim in this section is to give an interpretation of this equation in terms of Rankin-Cohen brackets. We have %
\[%
 \Phi_{1;2,1;12,0}(E_2,\Delta)\in\Delta\qm{4}{1}=\C\Delta E_4%
\]
hence %
\[%
 \Phi_{1;2,1;12,0}(E_2,\Delta)=\Delta E_4 %
\]
and
\[%
 \Phi_{1;4,0;12,0}(E_4,\Delta)\in\Delta\Modu{6}=\C\Delta E_6 %
\]
hence
\[%
 \Phi_{1;4,0;12,0}(E_4,\Delta)=4\Delta E_6 %
\]
so that %
\[%
 \Phi_{1;16,0;12,0}\left(\Phi_{1;2,1;12,0}(E_2,\Delta),\Delta\right)=\Delta\Phi_{1;4,0;12,0}(E_4,\Delta)=4\Delta^2E_6. %
\]
Next we compute %
\[%
 \Phi_{1;30,0;12,0}(\Delta^2E_6,\Delta)=\Delta^2\Phi_{1;6,0;12,0}(E_6,\Delta)\in\Delta^3\Modu{8}=\C\Delta^3E_4^2 %
\]
hence
\[%
 \Phi_{1;30,0;12,0}(\Delta^2E_6,\Delta)=6\Delta^3E_4^2=6\Delta\Phi_{1;2,1;12,0}(E_2,\Delta)^2 %
\]
and %
\[%
 \Phi_{1;30,0;12,0}\left(\Phi_{1;16,0;12,0}\left(\Phi_{1;2,1;12,0}(E_2,\Delta),\Delta\right),\Delta\right)=24\Delta\Phi_{1;2,1;12,0}(E_2,\Delta)^2. %
\]
This is \eqref{eq_chazyrc}. We deduce the usual form of the Chazy equation in the following way. From %
\[%
 K\coloneqq\Phi_{1;2,1;12,0}(E_2,\Delta)=E_2\De\Delta-12\De E_2\Delta=\Delta(E_2^2-12\De E_2) %
\]
we get %
\[%
 L\coloneqq\Phi_{1;16,0;12,0}(K,\Delta)=16K\Delta-12\De K\Delta=4\Delta^2(E_2^3-18E_2\De E_2+36\De^2E_2) %
\]
and since %
\begin{align*}
 \Phi_{1;30,0;12,0}(L,\Delta) &= 30L\De\Delta-12\De L\Delta\\%
&=24\Delta^3\left(E_2^4-24E_2^2\De E_2+72E_2\De^2E_2+36(\De E_2)^2-72\De^3E_2\right)%
\end{align*}
the equality $\Phi_{1;30,0;12,0}(L,\Delta)=24\Delta K^2$ gives the Chazy equation %
\[%
 2\De^3E_2-2E_2\De^2E_2+3(\De E_2)^2=0. %
\]

\bibliographystyle{amsalpha}
\bibliography{QuasimodRC}
\end{document}